\def\Bbb#1{{\fam\msbfam\relax#1}}
\font\fivemsb=msbm5
\font\sevenmsb=msbm7
\font\tenmsb=msbm10
\newfam\msbfam
\textfont\msbfam=\tenmsb
\scriptfont\msbfam\sevenmsb
\scriptscriptfont\msbfam\fivemsb

\def\spc{{\Bbb C}}
\def\spr{{\Bbb R}}
\def\vf{\varphi}
\def\wt{\widetilde}
\def\wh{\widehat}

\def\a{\alpha}
\def\b{\beta}

\def\l{\lambda}
\def\m{\mu}
\def\n{\nu}
\def\r{\rho}

\def\t{\theta}
\def\T{\Theta}
\def\x{\xi}
\def\y{\eta}
\def\ca{{\cal A}}

\def\ce{{\cal E}}

\def\cl{{\cal L}}

\def\la{\langle}
\def\ra{\rangle}
\def\lla{\langle\!\langle}
\def\rra{\rangle\!\rangle}

\documentstyle[11pt,amssymb,amscd,righttag]{amsart}
\theoremstyle{plain}
\newtheorem{theorem}{Theorem}[section]

\newtheorem{lemma}[theorem]{Lemma}

\theoremstyle{definition}
\newtheorem{definition}[theorem]{Definition}
\newtheorem{example}[theorem]{Example}

\numberwithin{equation}{section}

\begin{document}
%Topmatter
\title[]%
{Characterization of Hida Measures \\
in white noise analysis}
%\thanks{This paper is a considerably revised version 
%of an earlier preprint:
%``Characterization Theorem and Hida Measures in White 
%Noise Analysis" (1999).}

\maketitle
\begin{center}
{\sc Nobuhiro Asai
\footnote{Current address: 
International Institute for Advanced Studies,
Kizu, Kyoto, 619-0225, JAPAN.}}\\
{\it Graduate School of Mathematics\\
Nagoya University\\
Nagoya, 464-8602, JAPAN}
\end{center}
\smallskip
\begin{center}
{\sc Izumi Kubo}\\
{\it Department of Mathematics 
\\Graduate School of Science\\
Hiroshima University \\
Higashi-Hiroshima, 739-8526, JAPAN}
\end{center}
\begin{center}
and
\end{center}
\begin{center}
{\sc Hui-Hsiung Kuo}\\
{\it Department of Mathematics\\
Louisiana State University\\
Baton Rouge, LA 70803, USA}
\end{center}

\bigskip
\medskip
%%%%%%%%%%%%%%%%%%%%%%%%%%%%%%%%%%%%%%%%
% Introduction
%%%%%%%%%%%%%%%%%%%%%%%%%%%%%%%%%%%%%%%
\section{Introduction}
In the recent paper 
\cite{akk3} by Asai et al., the growth order of 
holomorphic functions on a nuclear space has been considered.
For this purpose,  certain classes of growth 
functions $u$ are introduced and many properties of 
Legendre transform of such functions are investigated. 
In \cite{akk4}, applying Legendre
transform of $u$ under the conditions (U0), (U2) and (U3) (see \S 2),
the Gel'fand triple
$$
[\ce]_u\subset (L^2)\subset [\ce]^{*}_u
$$ 
associated with a growth function $u$
is constructed.  

The main purpose of this work is to prove 
Theorem 4.4, so-called,  
the characterization theorem of Hida measures
(generalized measures).
As examples of such measures, we shall present
the Poisson noise measure and the Grey noise measure
in Example 4.5 and 4.6, respectively.

The present paper is organized as follows.
In \S 2, we give a quick review of 
some fundamental results in white noise 
analysis and introduce the notion of Legendre transform
utilized by Asai et al. in 
\cite{akk3},\cite{akk4}.
In \S 3, we simply cite some useful 
properties of the Legendre transform from 
\cite{akk3}.  In \S 4, 
we discuss the characterization of 
Hida measures (generalized measures).

\bigskip
\noindent
{\bf Acknowledgements.} 
N.~Asai wants to thank the Daiko
Foundation and the Kamiyama Foundation for research support.

%%%%%%%%%%%%%%%%%%%%%%%%%%%%%%%%%%%%%%%%
% Preliminaries
%%%%%%%%%%%%%%%%%%%%%%%%%%%%%%%%%%%%%%%%
\section{Preliminaries} \label{sec:1}
In this section, we will summarize well-known results in white noise 
analysis \cite{hkps},\cite{kuo96},\cite{ob} and notions from 
Asai et al.\cite{akk1},\cite{akk2},\cite{akk3},\cite{akk4}.
Complete details and further developments will be appeared in \cite{akk5}.
Some similar results have been obtained independently by Gannoun et al. 
\cite{ghor}.

Let $\ce_0$ be a real separale Hilbert space
with the norm $|\cdot|_0$.  Suppose
$\{|\cdot|_p\}_{p=0}^{\infty}$ is a sequence of 
densely defined inner product norms on $\ce_0$.
Let $\ce_p$ be the completion of $\ce$ with respect to
the norm $|\cdot|_p$.  In addition we assume 
\begin{itemize}
\item[(a)] There exists a constant $0<\r<1$ such that
$|\cdot|_{0}\leq \r|\cdot|_{1}\leq \cdots\leq
\r^{p}|\cdot|_{p} \leq \cdots$.
\item[(b)] For any $p\geq 0$, there exists $q\geq p$ such
that the inclusion $i_{q, p}: \ce_{q} \hookrightarrow
\ce_{p}$ is a Hilbert-Schmidt operator. 
\end{itemize}
Let $\ce'$ and $\ce_{p}'$ denote the dual spaces of $\ce$
and $\ce_{p}$, respectively. We can use the Riesz
representation theorem to identify $\ce_{0}$ with its
dual space $\ce_{0}'$. 
Let $\ce$ be the projective limit of 
$\{\ce_p ;\ p\geq 0\}$.
Then we get the following continuous
inclusions:
\begin{equation}
\ce \subset \ce_{p}  \subset \ce_{0}
  \subset \ce_{p}'  \subset \ce',
  \qquad p\geq 0.   \notag
\end{equation}
The above condition (b) says that $\ce$ is a nuclear
space and so $\ce \subset \ce_{0} \subset \ce'$ is a
Gel'fand triple. 

Let $\m$ be the standard Gaussian measure on $\ce'$
with the characteristic function given by
\begin{equation}
\int_{\ce'} e^{i\la x, \x\ra}\,d\m(x) =
  e^{-{1\over 2}|\x|_{0}^{2}}, \qquad \x\in\ce. \notag
\end{equation}
The probability space $(\ce', \m)$ is called 
a {\em white noise space} or {\em Gaussian space}. 
For simplicity, we will use
$(L^{2})$ to denote the Hilbert space of $\mu$-square
integrable functions on $\ce'$.
By the Wiener-It\^o theorem, each $\vf\in (L^{2})$ can be
uniquely expressed as 
\begin{equation} \label{eq:1-1}
\vf (x) = \sum_{n=0}^{\infty} I_{n}(f_{n})(x) =
\sum_{n=0}^{\infty} \la :\!x^{\otimes n}\!:, f_{n}\ra,
\qquad f_{n} \in \ce_{0}^{\wh\otimes n},
\end{equation}
where $I_{n}$ is the multiple Wiener integral of order $n$
and $:\!x^{\otimes n}\!:$ is the Wick tensor of $x\in\ce'$
(see \cite{kuo96}.) Moreover, the $(L^{2})$-norm
of $\vf$ is given by
\begin{equation} \label{eq:1-2}
\|\vf\|_{0} = \left(\sum_{n=0}^{\infty} n!|f_{n}|_{0}^{2}
\right)^{1/2}.
\end{equation}

Let $u\in C_{+,{1\over 2}}$ be the set of all positive 
continuous functions
on $[0,\infty)$ satisfying
$$
\lim_{r\to\infty}{\log u(r)\over\sqrt{r}}=\infty.
$$
In addition, we introduce conditions:
\begin{itemize}
\item[(U0)] 
$\inf_{r\geq 0} u(r) = 1$.
\smallskip
\item[(U1)] 
$u$ is increasing and $u(0)=1$.
\smallskip
\item[(U2)] 
$\lim_{r\to\infty} r^{-1} \log u(r) < \infty$.
\smallskip
\item[(U3)] 
$\log u(x^2)$ is convex on $[0,\infty)$.
\end{itemize}
Obviously, (U1) is a stronger condition than (U0). 

Let $C_{+, \log}$ denote the set of all positive continuous
functions $u$ on $[0, \infty)$ satisfying the condition:
\begin{equation}
\lim_{r\to\infty} {\log u(r) \over \log r}=\infty. \notag
\end{equation}
It is easy to see $C_{+,{1\over 2}}\subset C_{+,\log}$.

The {\em Legendre transform} $\,\ell_{u}$ of $u \in
C_{+, \log}$ is defined to be the function
\begin{equation}
\ell_{u}(t) = \inf_{r>0} {u(r) \over r^{t}}, \qquad
  t\in [0, \infty). \notag
\end{equation}
Some useful properties of the Legendre transform
will be refered in section \ref{sec:2}.

From now on, we take a function $u\in C_{+,{1\over 2}}$
satisfying (U0) (U2) (U3).  We shall constract 
a Gel'fand triple associated with $u$. 
For $\vf\in (L^{2})$ being represented by Equation
(\ref{eq:1-1}) and $p\geq 0$, define
\begin{equation} \label{eq:1-3}
\|\vf\|_{p, u} = \left(\sum_{n=0}^{\infty}{1\over \ell_u(n)}
   |f_{n}|_{p}^{2}\right)^{1/2}.
\end{equation}
Let $[\ce_{p}]_{u} = \{\vf\in (L^{2}); \,\|\vf\|_{p, u}
< \infty\}$. Define the space $[\ce]_{u}$ of {\em test
functions} on $\ce'$ to be the projective limit of
$\{[\ce_{p}]_{u}; \,p\geq 0\}$. The dual space
$[\ce]_{u}^{*}$ of $[\ce]_{u}$ is called the space of
{\em generalized functions} on $\ce'$.

Choose an appropriate $p_0$ such that 
$c\r^{2p_0}\sqrt{2}\leq 1$ for some $c$.  Then
two conditions (a) and (U2) imply 
that $[\ce_{p}]_{u} \subset (L^{2})$ 
for all $p\geq p_0$. Hence
$[\ce]_{u} \subset (L^{2})$ holds. By identifying
$(L^{2})$ with its dual space we get the following
continuous inclusions:
\begin{equation}
[\ce]_{u} \subset [\ce_{p}]_{u}
\subset (L^{2}) \subset
[\ce_{p}]_{u}^{*} \subset [\ce]_{u}^{*},
\qquad p \geq p_0,  \notag
\end{equation}
where $[\ce_{p}]_{u}^{*}$ is the dual space of
$[\ce_{p}]_{u}$. Moreover, $[\ce]_{u}$ is a nuclear
space and so $[\ce]_{u} \subset (L^{2}) \subset
[\ce]_{u}^{*}$ is a Gel'fand triple. Note that
$[\ce]_{u}^{*}=\cup_{p\geq 0} [\ce_{p}]_{u}^{*}$ and
for $p\geq p_0$,
$\,[\ce_{p}]_{u}^{*}$ is the completion of $(L^{2})$
with respect to the norm
\begin{equation} \label{eq:1-4}
\|\vf\|_{-p, (u)} = \left(\sum_{n=0}^{\infty} 
{(n!)^2\ell_u(n)}\,|f_{n}|_{-p}^{2}\right)^{1/2}.
\end{equation}

For $\x$ belonging to the complexification $\ce_{c}$ of
$\ce$, the renormalized exponential function $:\!e^{\la
\cdot, \x\ra}\!:$ is defined by
\begin{equation}
:\!e^{\la \cdot, \x\ra}\!: \, = \sum_{n=0}^{\infty}
{1 \over n!}\,\la :\!\cdot^{\otimes n}\!:,
\x^{\otimes n} \ra.  \notag
\end{equation}
Then we have the norm estimate,
\begin{equation}\label{eq:2-5}
\|:\!e^{\la \cdot, \x\ra}\!:\|^2_{-q,(u)}
=\sum_{n=0}^{\infty}
\ell_u(n)|\x|^{2n}_{-q}=:{\cal L}_u(|\x|^2_{-q}).
\end{equation}

For later uses, let us define the notion of 
{\it equivalent functions} here. 
\begin{definition}\label{def:equiv-func}
Two positive functions $f$ and $g$ on $[0, \infty)$ are
called {\em equivalent} if there exist constants
$c_{1}, c_{2}, a_{1}, a_{2} >0$ such that
\begin{equation}
c_{1}f(a_{1}r) \leq g(r) \leq c_{2}f(a_{2}r), \qquad
   \forall\> r\in [0, \infty).  \notag
\end{equation}
\end{definition}

%%%%%%%%%%%%%%%%%%%%%%%%%%%%%%%%%%%%%%
% Examples
%%%%%%%%%%%%%%%%%%%%%%%%%%%%%%%%%%%%%%
\begin{example}\label{example1}
\begin{equation} \label{eq:3-21}
g_{k}(r) = \exp\left[2\sqrt{r\log_{k-1}\sqrt{r}}\,\right],
\end{equation}
where $\log_k (r)$ is given by 
$$
\log_1(r)=\log(\max\{e,r\}), \quad
\log_k(r)=\log_1(\log_{k-1}(r)), \quad
k\geq 2.
$$
Then the function $g_{k}$ belongs to
$C_{+, 1/2}$ and satisfies conditions (U1) (U2) (U3). 
In the sense of Definition \ref{def:equiv-func},
the function $g_{k}$ is equivalent to the function
given by $$
u_k(r)=\sum_{n=0}^{\infty}
{1\over b_k(n)n!}r^n
$$
where $b_k(n)$ is the {\em k-th order Bell number}.
Hence we get the
Gel'fand triple,
$$
[\ce]_{g_{k}}\subset (L^{2}) \subset [\ce]_{g_{k}}^{*}
$$ 
known as the {\em CKS-space associated with $g_k$}, 
which is the same as the one defined by 
the k-th order Bell number $b_k(n)$.
See more details in 
\cite{akk1},\cite{akk2},\cite{akk3},
\cite{akk4},\cite{akk5},\cite{cks},\cite{kubo98},\cite{kks}.
\end{example}

\begin{example}\label{example2}
For $0\leq \b<1$, let $u$ be the function defined by
\begin{equation}
  u(r)= \exp\left[(1+\b)r^{1\over 1+\b}\right].  \notag
\end{equation}
It is easy to check that $u$ belongs to $C_{+, 1/2}$ and
satisfies conditions (U1) (U2) (U3). 
Hence this Gel'fand triple,
$$
(\ce)_{\b}\subset (L^{2}) \subset (\ce)_{\b}^{*}
$$
which is well-known as the {\em Hida-Kubo-Takenaka space} 
for $\b=0$ 
\cite{hkps},\cite{kon80},\cite{kta},
\cite{ktb},\cite{ob}
and the {\em Kondratiev-Streit space}
for a general $\b$ \cite{ks93}, \cite{kuo96}.
For $\b=1$ case, see 
\cite{kls96},\cite{ksw95},\cite{kswy98}.
\end{example}

\noindent
{\em Remark.}
We have the following chain of Gel'fand triples: 
$$
(\ce)_1\subset
[\ce]_{g_{k}}\subset [\ce]_{g_{l}}\subset 
(\ce)_{\b}\subset (\ce)_{\gamma}\subset (L^2)\subset
(\ce)^{*}_{\gamma}\subset (\ce)^{*}_{\b}
\subset [\ce]^{*}_{g_{l}}\subset [\ce]^{*}_{g_{k}}
\subset (\ce)^{*}_1
$$
where $0\leq \gamma\leq \b<1$ and $2\leq l\leq k$. 
\medskip
%%%%%%%%%%%%%%%%%%%%%%%%%%%%%%%%%%%%%%%%%%%%%%%%%%%%%%%%%%
% Section of Legendre transform
%%%%%%%%%%%%%%%%%%%%%%%%%%%%%%%%%%%%%%%%%%%%%%%%%%%%%%%%%
\section{Properties of Legendre transforms}
   \label{sec:2}

First we mention the following notions 
of concave and convex functions  
which will be used frequently. 
\begin{definition}
A positive function $f$ on $[0, \infty)$ is
called
\begin{itemize}
\item[(1)] {\em log-concave} if the function $\log f$ is
concave on $[0, \infty)$;
\item[(2)] {\em log-convex} if the function $\log f$ is
convex on $[0, \infty)$;
\item[(3)] {\em (log, exp)-convex} if the function
$\log f(e^{x})$ is convex on $\spr$;
\item[(4)] {\em (log, $x^{2}$)-convex} if the function
$\log f(x^{2})$ is convex on $[0, \infty)$.
\end{itemize}
\end{definition}

\smallskip
We will need the fact that if $f$ is log-concave,
then the sequence $\{f(n)\}_{n=0}^{\infty}$ is log-concave.
To check this fact, note that for any $t_{1}, t_{2}
\geq 0$ and $0\leq \l \leq 1$,
\begin{equation}
f(\l t_{1} + (1-\l)t_{2}) \geq f(t_{1})^{\l}\,
   f(t_{2})^{1-\l}.  \notag
\end{equation}
In particular, take $t_{1}=n, t_{2}=n+2$, and $\l=1/2$ to
get
\begin{equation}
f(n) f(n+2) \leq f(n+1)^{2}, \qquad
    \forall\> n\geq 0. \notag
\end{equation}
Hence the sequence $\{f(n)\}_{n=0}^{\infty}$ is
log-concave.

%Another fact we need later is from Proposition 2.3 (3) in
%\cite{akk3}: If $f$ is increasing and (log, $x^{2}$)-convex,
%then $f$ is (log, exp)-convex.

\smallskip
The next theorem is from Lemma 3.4 in \cite{akk3}.
\begin{theorem} \label{thm:2-1}
Let $u\in C_{+, \log}$. Then the Legendre transform
$\ell_{u}$ is log-concave. (Hence $\ell_{u}$ is
continuous on $[0, \infty)$ and the sequence
$\{\ell_{u}(n)\}_{n=0}^{\infty}$ is log-concave.)
\end{theorem}

From Theorem 2 (b) in \cite{akk1} we have the fact:
If $\{\a(n)/n!\}_{n=0}^{\infty}$ is log-concave and
$\a(0)=1$, then
\begin{equation}
\a(n+m) \leq {n+m \choose n} \a(n) \a(m), \qquad
\forall\, n, m\geq 0.  \notag
\end{equation}

By Theorem \ref{thm:2-1} the sequence $\{\ell_{u}(n)\}$
is log-concave. Hence we can apply the above fact to the
sequence $\a(n) = n!\ell_{u}(n)/\ell_{u}(0)$ to get the
next theorem.
\begin{theorem} \label{thm:2-2}
Let $u\in C_{+, \log}$. Then for all integers $n, m\geq 0$,
we have
\begin{equation}
\ell_{u}(0) \ell_{u}(n+m) \leq
    \ell_{u}(n) \ell_{u}(m).  \notag
\end{equation}
\end{theorem}

In the next theorem we state some properties of the
Legendre transform $\ell_{u}$ of a (log, exp)-convex
function $u$ in $C_{+, \log}$. It is from Lemmas 3.6 and
3.7 in \cite{akk3}.
\begin{theorem} \label{thm:2-3}
Let $u\in C_{+, \log}$ be (log, exp)-convex. Then
\begin{itemize}
\smallskip
\item[(1)] $\ell_{u}(t)$ is decreasing for large $t$,
\smallskip
\item[(2)] $\lim_{t\to\infty} \ell_{u}(t)^{1/t} =0$,
\smallskip
\item[(3)] $u(r) = \sup_{t\geq 0} \ell_{u}(t) r^{t}$
for all $r\geq 0$.
\end{itemize}
\end{theorem}

On the other hand, for a (log, $x^{2}$)-convex function
$u$ in $C_{+, \log}$, its Legendre transform $\ell_{u}$
has the properties in the next theorem from Lemmas 3.9
and 3.10 in \cite{akk3}. If in addition $u$ is increasing,
then $u$ is also (log, exp)-convex and hence $\ell_{u}$
has the properties in the above Theorem \ref{thm:2-3}.
\begin{theorem} \label{thm:2-4}
Let $u\in C_{+, \log}$. We have the assertions:
\begin{itemize}
\smallskip
\item[(1)] $u$ is (log, $x^{2}$)-convex if and only if
$\ell_{u}(t) t^{2t}$ is log-convex.
\smallskip
\item[(2)] If $u$ is (log, $x^{2}$)-convex, then for any
integers $n, m\geq 0$,
\begin{equation}
\ell_{u}(n) \ell_{u}(m) \leq \ell_{u}(0) 2^{2(n+m)}
  \ell_{u}(n+m).  \notag
\end{equation}
\end{itemize}
\end{theorem}

Now, suppose $u\in C_{+, \log}$ and assume that
$\lim_{n\to\infty} \ell_{u}(n)^{1/n} =0$. We define the
$L$-{\em function} $\cl_{u}$ of $u$ by
\begin{equation} \label{eq:2-1}
  \cl_{u}(r) = \sum_{n=0}^{\infty} \ell_{u}(n) r^{n}.
\end{equation}

Note that $\cl_{u}$ is an entire function. By Theorem
\ref{thm:2-3} (2), $\ell_{u}$ is defined for any
(log, exp)-convex function $u$ in $\in C_{+, \log}$.
Moreover, we have the next theorem from Theorem 3.13 in
\cite{akk3}.
\begin{theorem} \label{thm:2-5}
(1) Let $u\in C_{+, \log}$ be (log, exp)-convex. Then
its $L$-function $\cl_{u}$ is also (log, exp)-convex and
for any $a>1$,
\begin{equation}
\cl_{u}(r) \leq {ea \over \log a} u(ar), \qquad
\forall r\geq 0.  \notag
\end{equation}
\par\noindent
(2) Let $u\in C_{+, \log}$ be increasing and
(log, $x^{2}$)-convex. Then there exists a constant $C$
such that
\begin{equation}
u(r) \leq C \cl_{u}(2^{2} r), \qquad \forall r\geq 0.
    \notag
\end{equation}
\end{theorem}

Recall from Proposition 2.3 (3) in \cite{akk3}:
If $f$ is increasing and (log, $x^{2}$)-convex for
some $k>0$, then $f$ is (log, exp)-convex. Hence the above
Theorem \ref{thm:2-5} yields the next theorem.
\begin{theorem} \label{thm:2-a}
Let $u\in C_{+, \log}$ be increasing and
(log, $x^{2}$)-convex. Then the functions
$u$ and $\cl_{u}$ are equivalent.
\end{theorem}

In the next section \ref{sec:4}, 
we will consider the characterization of Hida measures
(generalized measures).
We prepare two lemmas for this purpose.
The proof of Lemma \ref{lem:4-1} is simple application of 
Theorem \ref{thm:2-4} so that we just state it without proof.
\begin{lemma} \label{lem:4-1}
Suppose $u\in C_{+, \log}$ is (log, $x^{2}$)-convex. Then
\begin{equation}
\cl_{u}(r)^{2} \leq \ell_{u}(0) \cl_{u}\big(2^{3}r\big),
  \qquad \forall r\in [0, \infty).
\end{equation}
\end{lemma}

\noindent
{\em Remark.} Note that $\cl_{u}(r)\geq \ell_{u}(0)$ for
all $r\geq 0$. Hence we have
\begin{equation}
\ell_{u}(0) \cl_{u}(r) \leq \cl_{u}(r)^{2} \leq \ell_{u}(0)
\cl_{u}\big(2^{3}r\big), \qquad \forall r\in [0, \infty).
\notag
\end{equation}
Thus $\cl_{u}$ and $\cl_{u}^{2}$ are equivalent for any
(log, $x^{2}$)-convex function $u\in C_{+, \log}$. If, in
addition, $u$ is increasing, then $u$ and $\cl_{u}$ are
equivalent by Theorem \ref{thm:2-a}. It follows that $u$ and
$u^{2}$ are equivalent for such a function $u$.

The next Lemma \ref{lem:4-2} can be obtained from 
Theorem \ref{lem:4-1} and Lemma \ref{thm:2-5}.
\begin{lemma} \label{lem:4-2}
Suppose $u\in C_{+, \log}$ is increasing and
(log, $x^{2}$)-convex. Then for any $a>1$, we have
\begin{equation} \label{eq:4-11}
\cl_{u}(r) \leq \sqrt{\ell_{u}(0){ea \over \log a}}\>
   u\big(a2^{3}r\big)^{1/2}.
\end{equation}
\end{lemma}

%%%%%%%%%%%%%%%%%%%%%%%%%%%%%%%%%%%%%%%%%%%%%%%%%%%%%%%%%
% Section of Hida Measure
%%%%%%%%%%%%%%%%%%%%%%%%%%%%%%%%%%%%%%%%%%%%%%%%%%%%%%%%% 
\section{Characterization of Hida Measures} \label{sec:4}

Before going to the main theorem, we need to 
%In this section we will 
introduce another equivalent
family of norms on $[\ce]_{u}$, i.e.,
$\{\|\cdot\|_{\ca_{p, u}};\,p\geq 0\}$. This family of
norms is intrinsic in the sense that $\|\vf\|_{\ca_{p, u}}$
is defined directly in terms of the analyticity and growth
condition of $\vf$.

\smallskip
First, it is well-known that each test function $\vf$ in
$[\ce]_{u}$ has a unique analytic extension (see \S 6.3
of \cite{kuo96}) given by
\begin{equation} \label{eq:4-1}
\vf(x) = \lla :\!e^{\la\cdot, x\ra}\!:, \T\vf\rra,
\qquad x\in \ce_{c}',
\end{equation}
where $\T$ is the unique linear operator taking
$e^{\la\cdot, \x\ra}$ into $:\!e^{\la\cdot, \x\ra}\!:$ for
all $\x\in \ce_{c}$. 
%This operator turns out to be the
%same as $\cg_{i, 1}$ defined in Equation (\ref{eq:fg})
%in section section \ref{sec:5}. 
By Theorem 6.2 in \cite{kuo96} with minor modifications,
$\T$ is shown to be a continuous linear operator
from $[\ce]_{u}$ into itself.
Note that we still assume conditions (U0), (U2) and (U3)
on $u$ given in section \ref{sec:1}. 

Now, let $p\geq 0$ be any fixed number. Choose $p_{1}>p$
such that $2\r^{2(p_{1}-p)}\leq 1$. Then use Equations
\eqref{eq:4-1}, \eqref{eq:2-5} 
and Theorem \ref{thm:2-5} to get
\begin{equation}
|\vf(x)| \leq \|\T\vf\|_{p_{1}, u}\,\|\!:\!e^{\la\cdot,
  x\ra}\!:\!\|_{-p_{1}, (u)}
  \leq \|\T\vf\|_{p_{1}, u}\, \sqrt{2e \over \log 2}\,
   u\big(2|x|_{-p_{1}}^{2}\big)^{1/2}.  \notag
\end{equation}
Note that $2|x|_{-p_{1}}^{2}\leq 2\r^{2(p_{1}-p)}
|x|_{-p}^{2} \leq |x|_{-p}^{2}$ by the above choice of
$p_{1}$. Since $u$ is an increasing function, we see that
\begin{equation}
  |\vf(x)| \leq \|\T\vf\|_{p_{1}, u}\, \sqrt{2e \over
  \log 2}\,u\big(|x|_{-p}^{2}\big)^{1/2}.  \notag
\end{equation}

But $\T$ is a continuous linear operator from $[\ce]_{u}$
into itself. Hence there exist positive constants $q$
and $K_{p, q}$ such that $\|\T\vf\|_{p_{1}, u}\leq
K_{p, q}\|\vf\|_{q, u}$. Therefore,
\begin{equation} \label{eq:4-2}
  |\vf(x)| \leq C_{p, q} \|\vf\|_{q, u}\,u\big(|x|_{-p}^{2}
  \big)^{1/2},  \qquad x \in \ce_{p, c}',
\end{equation}
where $C_{p, q}=K_{p, q}\sqrt{2e/\log 2}$. This is the
growth condition for test functions.

Being motivated by Equation (\ref{eq:4-2}), we define
\begin{equation} \label{eq:4-3}
\|\vf\|_{\ca_{p, u}} = \sup_{x\in\ce_{p, c}'} |\vf(x)|
\,u\big(|x|_{-p}^{2}\big)^{-1/2}.
\end{equation}
Obviously, $\|\cdot\|_{\ca_{p, u}}$ is a norm on
$[\ce]_{u}$ for each $p\geq 0$.

\begin{theorem} \label{thm:4-1}
Suppose $u\in C_{+, 1/2}$ satisfies conditions (U1) (U2)
(U3). Then the families of norms $\{\|\cdot\|_{\ca_{p, u}};
\,p\geq 0\}$ and $\{\|\cdot\|_{p, u};\,p\geq 0\}$ are
equivalent, i.e., they generate the same topology on
$[\ce]_{u}$.
\end{theorem}

\noindent
{\em Remark.} This theorem 
gives an alternative construction of test
functions. This idea is due to Lee \cite{lee}, see also
\S 15.2 of \cite{kuo96}. For $p\geq 0$, let
$\ca_{p, u}$ consist of all functions $\vf$ on $\ce_{c}'$
satisfying the conditions:
\begin{itemize}
\item[(a)] $\vf$ is an analytic function on $\ce_{p, c}'$.
\item[(b)] There exists a constant $C\geq 0$ such that
\begin{equation}
|\vf(x)| \leq C u\big(|x|_{-p}^{2}\big)^{1/2}, \qquad
  \forall x \in \ce_{p, c}'.  \notag
\end{equation}
\end{itemize}

For each $\vf\in \ca_{p, u}$, define $\|\vf\|_{\ca_{p, u}}$
by Equation (\ref{eq:4-3}). Then $\ca_{p, u}$ is a Banach
space with norm $\|\cdot\|_{\ca_{p, u}}$. Let $\ca_{u}$ be
the projective limit of $\{\ca_{p, u};\,p\geq 0\}$. We can
use the above theorem to conclude that $\ca_{u}=[\ce]_{u}$
as vector spaces with the same topology. Here the equality
$\ca_{u}=[\ce]_{u}$ requires the use of analytic extensions
of test functions in $[\ce]_{u}$, which exists in view of
Equation (\ref{eq:4-1}).

\begin{pf}
Let $p\geq 0$ be any given number. We have already shown
that there exist constants $q>p$ and $C_{p, q}\geq 0$
such that Equation (\ref{eq:4-2}) holds. It follows that
\begin{equation}
\|\vf\|_{\ca_{p, u}} = \sup_{x\in\ce_{p, c}'} |\vf(x)|
\,u\big(|x|_{-p}^{2}\big)^{-1/2}
  \leq C_{p, q} \|\vf\|_{q, u}.  \notag
\end{equation}
Hence for any $p\geq 0$, there exist constants $q>p$ and
$C_{p, q}\geq 0$ such that
\begin{equation} \label{eq:4-4}
  \|\vf\|_{\ca_{p, u}} \leq C_{p, q} \|\vf\|_{q, u},
  \qquad \forall \vf\in [\ce]_{u}.
\end{equation}

To show the converse, first note that by condition (U2)
there exist constants $c_{1}, c_{2}>0$ such that
$u(r)\leq c_{1}e^{c_{2}r}, \> r\geq 0$. Next note that
by Fernique's theorem 
(see \cite{fer}, \cite{kuo75}, \cite{kuo96}) we have
\begin{equation}
\int_{\ce'} e^{2c_{2}|x|_{-\l}^{2}}\,d\m(x) < \infty
\qquad \text{for all large~} \l.  \notag
\end{equation}

Now, let $p\geq 0$ be any given number. Choose $q>p$
large enough such that
\begin{equation} \label{eq:4-5}
4e^{2}\|i_{q, p}\|_{HS}^{2}<1, \qquad
\int_{\ce'} e^{2c_{2}|x|_{-q}^{2}}\,d\m(x) < \infty.
\end{equation}
With this choice of $q$ we will show below that
\begin{equation} \label{eq:4-6}
\|\vf\|_{p, u} \leq L_{p, q} \|\vf\|_{\ca_{q, u}},
\qquad \forall \vf\in [\ce]_{u},
\end{equation}
where $L_{p, q}$ is the constant given by
\begin{equation} \label{eq:4-a}
L_{p, q} = \sqrt{c_{1}}\left(1-4e^{2}
  \|i_{q, p}\|_{HS}^{2}\right)^{-1/2} \int_{\ce'}
  e^{2c_{2}|x|_{-q}^{2}}\,d\m(x).
\end{equation}
Observe that the theorem follows from Equations
(\ref{eq:4-4}) and (\ref{eq:4-6}).

Finally, we prove Equation (\ref{eq:4-6}). 
Let $\vf\in [\ce]_{u}$.
% and $F={\cal S}\vf$. 
Then we can use an integral form of S-transform (see \cite{kuo96})
given by 
\begin{equation}
F(\x) = S\vf(\x) = \int_{\ce'} \vf(x+\x)\,d\m(x), \qquad
  \x\in\ce_{c}.  \notag
\end{equation}
Hence for the above choice of $q$, we have
\begin{align}
|F(\x)| & \leq \int_{\ce'} |\vf(x+\x)|\,d\m(x)
      \notag  \\
  & \leq \int_{\ce'} \left(|\vf(x+\x)|\,
     u\big(|x+\x|_{-q}^{2}\big)^{-1/2}\right)
    u\big(|x+\x|_{-q}^{2}\big)^{1/2}\,d\m(x)
      \notag  \\
  & \leq \|\vf\|_{\ca_{q, u}} \int_{\ce'}
    u\big(|x+\x|_{-q}^{2}\big)^{1/2}\,d\m(x).
      \notag
\end{align}
Here by condition (U1), we have $u(r)^{1/2}\leq
u(r)$ for all $r\geq 0$. Therefore,
\begin{equation} \label{eq:4-7}
|F(\x)| \leq \|\vf\|_{\ca_{q, u}} \int_{\ce'}
    u\big(|x+\x|_{-q}^{2}\big)\,d\m(x).
\end{equation}
By condition (U3), we have
\begin{equation}
u\left(\big({\textstyle{1\over 2}}r_{1}+{\textstyle{1\over
     2}}r_{2}\big)^{2}\right)
\leq u\big(r_{1}^{2}\big)^{1/2}\,
   u\big(r_{2}^{2}\big)^{1/2},
\qquad \forall r_{1}, r_{2} \geq 0. \notag
\end{equation}
Put $r_{1}=2|x|_{-q}$ and $r_{2}=2|\x|_{-q}$ to get
\begin{align}
u\big(|x+\x|_{-q}^{2}\big)
& \leq u\left(\big({\textstyle{1\over 2}} 2|x|_{-q}
   +{\textstyle{1\over 2}} 2|\x|_{-q}\big)^{2}
       \right)      \notag \\
& \leq u\big(4|x|_{-q}^{2}\big)^{1/2}\,
    u\big(4|\x|_{-q}^{2}\big)^{1/2}.   \notag
\end{align}
Then integrate over $\ce'$ to obtain the inequality:
\begin{equation} \label{eq:4-8}
\int_{\ce'} u\big(|x+\x|_{-q}^{2}\big)\,d\m(x)
\leq u\big(4|\x|_{-q}^{2}\big)^{1/2} \int_{\ce'}
  u\big(4|x|_{-q}^{2}\big)^{1/2}\,d\m(x).
\end{equation}
Put Equation (\ref{eq:4-8}) into Equation (\ref{eq:4-7})
to get
\begin{equation} \label{eq:4-9}
|F(\x)| \leq \|\vf\|_{\ca_{q, u}}
  u\big(4|\x|_{-q}^{2}\big)^{1/2} \int_{\ce'}
  u\big(4|x|_{-q}^{2}\big)^{1/2}\,d\m(x).
\end{equation}

Now, by the inequality $u(r)\leq c_{1}e^{c_{2}r}$, we have
\begin{equation} \label{eq:4-10}
\int_{\ce'} u\big(4|x|_{-q}^{2}\big)^{1/2}\,d\m(x)
   \leq  \sqrt{c_{1}} \int_{\ce'} e^{2c_{2}|x|_{-q}^{2}}
    \,d\m(x),
\end{equation}
which is finite by the choice of $q$ in Equation
(\ref{eq:4-5}).

From Equations (\ref{eq:4-9}) and (\ref{eq:4-10}),
we see that
\begin{equation}
  |F(\x)| \leq \|\vf\|_{\ca_{q, u}} \sqrt{c_{1}}
  \left(\int_{\ce'} e^{2c_{2}|x|_{-q}^{2}}\,d\m(x)\right)
    u\big(4|\x|_{-q}^{2}\big)^{1/2},
    \qquad \x\in\ce_{c}. \notag
\end{equation}
With this inequality and the choice of $q$ in Equation
(\ref{eq:4-5}) we can apply Lemma \ref{lem:3-2} (see below)
and Equation \eqref{eq:1-3} to
show that for any $\vf\in [\ce]_{u}$,
\begin{equation}
\|\vf\|_{q, u} \leq L_{p, q} \|\vf\|_{\ca_{q, u}}, \notag
\end{equation}
where $L_{p, q}$ is given by Equation(\ref{eq:4-a}). Thus
the inequality in Equation (\ref{eq:4-6}) holds and so the
proof is completed.
\end{pf}

In the proof of the prevous theorem, 
we have used the next lemma from \cite{akk3}.
\begin{lemma}[\cite{akk3}] \label{lem:3-2}
Suppose $u\in C_{+, 1/2}$ satisfies conditions (U1) (U2)
(U3). Let $F$ be a complex-valued function on $\ce_{c}$
satisfying the conditions:
\begin{itemize}
\item[(1)] For any $\x, \y \in \ce_{c}$, the function
$F(z\x+\y)$ is an entire function of $z\in\spc$.
\item[(2)] There exist constants $K, a, p\geq 0$ such
that
\begin{equation}
|F(\x)| \leq K u\big(a|\x|_{-p}^{2}\big)^{1/2},
  \qquad \x\in\ce_{c}.  \notag
\end{equation}
\end{itemize}
Let $q\in [0, p)$ be a number such that
$ae^{2}\|i_{p, q}\|_{HS}^{2} < 1$. 
Then there exist functions $f_n\in\ce^{\widehat{\otimes}}_{q,\spc}$ 
such that $F(\x)=\sum_{n=0}^{\infty}\la f_n,\x^{\widehat{\otimes}}\ra$
and 
\begin{equation} 
|f_n|^2_q\leq K(ae^2\|i_{p,q}\|^2_{HS})^n\ell_u(n).
\end{equation}
\end{lemma}

\begin{definition}
A measure $\n$ on $\ce'$ is called a {\em Hida measure}
associated with $u$ if $[\ce]_{u}\subset L^{1}(\n)$ and
the linear functional $\vf\mapsto \int_{\ce'}
\vf(x)\,d\n(x)$ is continuous on $[\ce]_{u}$. 
\end{definition}

In this case,
$\n$ induces a generalized function, denoted by $\wt\n$,
in $[\ce]_{u}^{*}$ such that
\begin{equation} \label{eq:4-12}
\lla \wt\n, \vf\rra = \int_{\ce'} \vf(x)\,d\n(x),
  \qquad \vf \in [\ce]_{u}.
\end{equation}

\begin{theorem} \label{thm:4-2}
Suppose $u\in C_{+, 1/2}$ satisfies conditions (U1) (U2)
(U3). Then a measure $\n$ on $\ce'$ is a Hida measure with
$\wt\n\in [\ce]_{u}^{*}$ if and only if $\n$ is supported
by $\ce_{p}'$ for some $p\geq 0$ and
\begin{equation} \label{eq:4-13}
\int_{\ce_{p}'} u\big(|x|_{-p}^{2}\big)^{1/2}\,d\n(x)
  < \infty.
\end{equation}
\end{theorem}

\noindent
{\em Remarks.} (a) The integrability condition in the
theorem can be replaced by
\begin{equation}
\int_{\ce_{p}'} u\big(|x|_{-p}^{2}\big)\,d\n(x)
  < \infty. \notag
\end{equation}
To verify this fact, just note that $u$ and $u^{2}$ are
equivalent (from the Remark of Lemma \ref{lem:4-1}) and
$|x|_{-q}\leq \r^{q-p} |x|_{-p}$ for $0\leq p\leq q$
and $x\in \ce_{p}'$.

\smallskip
(b) This theorem is due to Lee \cite{lee} for the case
$u(r)=e^{r}$. See \S 15.2 of the book \cite{kuo96} for the
case $u(r)=\exp\big[(1+\b) r^{1\over 1+\b}\big],
\,0\leq \b<1$.  In the case of $\b=1$, we need 
special treatment since our Legendre transform 
method should be modified.  In 
order to take care of $\b=1$ case,
we have to remove the assumption 
$$
	\lim_{r\to\infty}{\log u(r)\over \sqrt{r}}=\infty
$$
on $u$ introduced in \S 2, for example. 
It will be discussed in the future.
On the other hand, there are references 
\cite{ksw95},\cite{kswy98} discussed 
this case with a diffrent way from our point of view.
\begin{pf}
To prove the sufficiency, suppose $\n$ is supported by
$\ce_{p}'$ for some $p\geq 0$ and Equation (\ref{eq:4-13})
holds. Then for any $\vf\in [\ce]_{u}$,
\begin{align}
\int_{\ce'} |\vf(x)|\,d\n(x)
  & = \int_{\ce_{p}'} |\vf(x)|\,d\n(x)  \notag  \\
  & = \int_{\ce_{p}'} \left(|\vf(x)| u\big(|x|_{-p}^{2}
      \big)^{-1/2}\right) u\big(|x|_{-p}^{2}\big)^{1/2}
       \,d\n(x)  \notag  \\
  & \leq \|\vf\|_{\ca_{p, u}} \int_{\ce_{p}'}
   u\big(|x|_{-p}^{2}\big)^{1/2}\,d\n(x). \label{eq:4-14}
\end{align}

By Theorem \ref{thm:4-1}, $\{\|\cdot\|_{\ca_{p, u}};
\,p\geq 0\}$ and $\{\|\cdot\|_{p, u};\,p\geq 0\}$ are
equivalent. Hence Equation (\ref{eq:4-14}) implies that
$[\ce]_{u}\subset L^{1}(\n)$ and the linear functional
\begin{equation}
  \vf\longmapsto \int_{\ce'} \vf(x)\,d\n(x), \qquad
  \vf\in [\ce]_{u}, \notag
\end{equation}
is continuous on $[\ce]_{u}$. Thus $\n$ is a Hida measure
with $\wt\n$ in $[\ce]_{u}^{*}$.

To prove the necessity, suppose $\n$ is a Hida measure
inducing a generalized function $\wt\n \in [\ce]_{u}^{*}$.
Then for all $\vf \in [\ce]_{u}$,
\begin{equation} \label{eq:4-15}
\lla \wt\n, \vf \rra = \int_{\ce'} \vf(x)\,d\n(x).
\end{equation}
Since $\{\|\cdot\|_{\ca_{p, u}}; \,p\geq 0\}$ and
$\{\|\cdot\|_{p, u};\,p\geq 0\}$ are equivalent, the
linear functional $\vf\mapsto \lla \wt\n, \vf \rra$ is
continuous with respect to $\{\|\cdot\|_{\ca_{p, u}};
\,p\geq 0\}$. Hence there exist constants $K, q\geq 0$
such that for all $\vf \in [\ce]_{u}$,
\begin{equation} \label{eq:4-16}
\big| \lla \wt\n, \vf \rra\big| \leq
   K\|\vf\|_{\ca_{q, u}}.
\end{equation}
Note that by continuity, Equations (\ref{eq:4-15}) and
(\ref{eq:4-16}) also hold for all $\vf\in \ca_{q, u}$
defined in the Remark of Theorem \ref{thm:4-1}.

Now, with this $q$, we define a function $\t$ on
$\ce_{q, c}'$ by
\begin{equation}
\t (x) = \cl_{u}\big(2^{-4} \la x, x\ra_{-q}\big),
\qquad x\in \ce_{q, c}',  \notag
\end{equation}
where $\la \cdot, \cdot\ra_{-q}$ is the bilinear pairing
on $\ce_{q, c}'$. Obviously, $\t$ is analytic on
$\ce_{q, c}'$. On the other hand, apply Lemma \ref{lem:4-2}
with $a=k=2$ to get
\begin{equation}
|\t(x)| \leq \cl_{u}\big(2^{-4}|x|_{-q}^{2}\big) \leq
  \sqrt{2e \over \log 2}\> u\big(|x|_{-q}^{2}\big)^{1/2},
   \qquad \forall x\in\ce_{q, c}'.  \notag
\end{equation}
This shows that $\t\in \ca_{q, u}$ and we have
\begin{equation}
\|\t\|_{\ca_{q, u}} \leq \sqrt{2e \over \log 2}. \notag
\end{equation}
Then apply Equation (\ref{eq:4-16}) to the function $\t$,
\begin{equation}
  \big|\lla \wt\n, \t \rra\big| \leq K\|\t\|_{\ca_{q, u}}
  \leq K \sqrt{2e \over \log 2}. \notag
\end{equation}
Therefore, from Equation (\ref{eq:4-15}) with $\vf=\t$
we conclude that
\begin{equation} \label{eq:4-17}
\left|\int_{\ce'} \t(x)\,d\n(x)\right| \leq
   K \sqrt{2e \over \log 2}.
\end{equation}

Note that $\t(x)=\cl_{u}\big(2^{-4}|x|_{-q}^{2}\big)$ for
$x\in\ce'$. Hence Equation (\ref{eq:4-17}) implies that
\begin{equation}
\int_{\ce'} \cl_{u}\big(2^{-4}|x|_{-q}^{2}\big)\,d\n(x)
  < \infty.  \notag
\end{equation}
But $u(r) \leq C \cl_{u}(4r)$ from Theorem \ref{thm:2-5}
(2) with $k=2$. Therefore,
\begin{equation}
  \int_{\ce'} u\big(2^{-6}|x|_{-q}^{2}\big)\,d\n(x)
  < \infty.  \notag
\end{equation}
Now, choose $p>q$ large enough such that $\r^{2(p-q)}
\leq 2^{-6}$. Then $|x|_{-p}^{2} \leq 2^{-6}|x|_{-q}^{2}$.
Recall that $u$ is increasing. Hence
\begin{equation}
  \int_{\ce'} u\big(|x|_{-p}^{2}\big)\,d\n(x)
  < \infty.  \notag
\end{equation}
Note that $u(r)\geq 1$ and so $u(r)^{1/2}(r) \leq u(r)$.
Thus we conclude that
\begin{equation}
\int_{\ce'} u\big(|x|_{-p}^{2}\big)^{1/2}\,d\n(x)
  < \infty.  \notag
\end{equation}
This inequality implies that $\n$ is supported by
$\ce_{p}'$ and Equation (\ref{eq:4-13}) holds.
\end{pf}

%%%%%%%%%%%%%%%%%%%%%%%%%%%%%%%%%%%%%%%%%%%%%%%
% Examples of Hida measures
%%%%%%%%%%%%%%%%%%%%%%%%%%%%%%%%%%%%%%%%%%%%%%%
\begin{example}\label{example-poisson}
(Poisson noise measure)\\
Let ${\cal P}$ be the Poisson measure on $\ce^{*}$ 
given by 
\begin{equation*}
\exp\Bigl(\int_{\spr}(e^{i\x(t)}-1)dt\Bigr)
=\int_{\ce^{*}}e^{i\la x,\x\ra}{\cal P}(dx),
\quad \x\in\ce^{*}.
\end{equation*}
It has been shown \cite{cks} that the Poisson noise 
measure induces a generalized function in $[\ce]^{*}_{g_2}$.
%the CKS-space associated with $g_2$. 
Thus by Theorem \ref{thm:4-2} and Example \ref{example1}
we have the 
integrability condition
\begin{equation*}
\int_{\ce^{*}_p}
\exp\Bigl(|x|_{-p}\sqrt{\log |x|_{-p}}\Bigr)
{\cal P}(dx)<\infty
\end{equation*}
for some $p$. 
\end{example}

\begin{example}\label{example-grey}
(Grey noise measure)\\
Let $0<\lambda\leq 1$.
The grey noise measure on $\ce^{*}$ 
is the measure $\nu_{\lambda}$ having
the characteristic function
\begin{equation*}
L_{\lambda}(|\x|^2_0)
=\int_{\ce^{*}}e^{i\la x,\x\ra}\nu_{\lambda}(dx),
\quad \x\in \ce,
\end{equation*}
where $L_{\lambda}(t)$ is the Mittag-Leffler
function with parameter $\lambda$;
\begin{equation*}
L_{\lambda}(t)=\sum_{n=0}^{\infty}
{(-t)^n\over\Gamma(1+\lambda n)}.
\end{equation*}
Here $\Gamma$ is the Gamma function.
This measure was introduced by Schneider \cite{schneider}.
It is shown in \cite{kuo96} that
$\nu_{\lambda}$ is a Hida measure which 
induces a generalized function 
$\varPhi_{\nu_{\lambda}}$ in $(\ce)^{*}_{1-\lambda}$.
Therefore by Theorem \ref{thm:4-2} and 
Example \ref{example2}
the grey noise measure $\nu_{\lambda}$ 
satisfies 
\begin{equation*}
\int_{\ce^{*}_{p}}
\exp\bigl({1\over 2}(2-\lambda)
|x|_{-p}^{{2\over 2-\lambda}}\bigr)\nu_{\lambda}(dx)
<\infty
\end{equation*}
for some $p$. 
\end{example}

%%%%%%%%%%%%%%%%%%%%%%%%%%%%%%%%%%%%%%%%%%%%%%%%%%%%%
% References
%%%%%%%%%%%%%%%%%%%%%%%%%%%%%%%%%%%%%%%%%%%%%%%%%%%%

\end{document}